\input amstex
\magnification=1440
\hcorrection{-0,5cm}

\documentstyle{amsppt}
\pageheight{37pc}
\pagewidth{28,5pc}
\define\1{\hbox{\rm 1}\!\! @,@,@,@,@,\hbox{\rm I}}

\vskip 0,5cm

\topmatter
\title
{
On invariance of domains with smooth boundaries
with respect to stochastic differential equations
}
\endtitle
\footnotetext{This research was supported (in part) by the
Ministry of Education and Science of Ukraine, project
No 01.07/103.}
\author
{Vitalii A. Gasanenko}
\endauthor
\address
{Institute of Mathematics,National Academy of
Science of Ukraine,Tereshchenkivska 3, 252601, Kiev, Ukraine}
\endaddress
\endkeywords
\subjclass
{60 J 60}
\endsubjclass
\email
{gs\@imath.kiev.ua}
\endemail
\abstract
We prove constructible sufficient conditions of lack of exit
by$\quad$ solutions of stochastic differential Ito's equations
from domains with smooth boundaries.
\endabstract
\endtopmatter
\rightheadtext{
On invariance of domaines with smooth boundaries
}
\document

Consider a stochastic differential equation for process
$\xi(t)\in R^{n}$.

$$
d\xi(t)=a(t,\xi)dt+\sum\limits_{k=1}^{n}b_{k}(t,\xi)dw_{k}
,\quad \xi(0)=\xi_{0}\eqno (1)
$$

here

$a(t,x):=(a_{i}(t,x), \, 1\leq i\leq n),\quad b_{k}(t,x):=(b_{ki}(x),\,
1\leq i\leq n),\, x\in R^{n}$.

It is assumed that there is such constant
~$L$ that for functions ~$b_{ij}(t,x),~ a_{i}(t,x)$ 

\noindent the following
conditions take place

$$
|a(s,x)- a(s,y)| +
\sum\limits_{k}^{n}|b_{k}(s,x)- b_{k}(s,y)|
\leq L|x-y|,
$$

$$
|a(s,x)|^{2} + \sum\limits_{k}^{n}|b_{k}(s,x)|^{2}
\leq L^{2}(1+|x|^{2}).\eqno(2)
$$

for all $x,y\in R^{n}$.

Here~ $|\cdot|$,~ be norm (length) of vector.

It follows from  [ 7, p.480] that this is sufficient conditions for
existence of unique solution of (1).

   Let there be given a measurable set  $K\in R^{n}$. A set $K$ is said
to be invariant set of equation  (1) , if under condition
~$P(\xi_{0}\in K)=1$~the following equalities hold

$$
P(\xi(t)\in K)=1\quad \hbox{for all}\quad
 t\geq 0.
\eqno(3)
$$

The property (3) of trajectories of solution of (1) sometimes is called
viability.
The necessary and sufficient conditions of viability
was proved for the first time in [1].
Such conditions was proved for more general constructions of equations
in [2]. The methods of investigations of these articles are different but
the set ~$K$~is the same: it is  convex and closure. 

This problem was reduced to viability
of ordinary differential equations with help approximation theorems
Ikeda- Nakao- Yamato for homogeneous stochastic differential equations in [4].

The conditions of viability were formulated
in terms of asymptotic behavior of  distance to considered closed set.
The analogy conditions of viability were proved for inhomogeneous
stochastic differential equations and relative closed sets in [5].
We observe that test of conditions in terms of distance to sets
requires the  additional investigations. They are checked for convex
sets 

effectively. For example, it was done in [5].

 It was proved necessary and sufficient conditions or only sufficient

conditions of viadility (3) in [3,6] by probabilistic methods for the
specific domains ~$K$.

Our purpose is to obtain verifiable sufficient
conditions of viability (3) for domains with smooth boundaries.
Our method of investigation be different from other. It is based
on the use of Ostrogradskii- Gauss theorem.

Consider a closed set $K$ in $R^{n}$ with boundary $S$( or $\partial K$).

Let  $U(x,r)$ denote open ball with center in point $x$ and with radius
$r$. The union of balls with centers in $K$ is called $\epsilon$ - neighborhood
 $K_{\epsilon}$ of the set $K$ :
$K_{\epsilon}=\bigcup\limits_{x\in K}U(x,\epsilon)$.
We will denote by ~~ $S_{\epsilon}$ the boundary of ~$K_{\epsilon}$.

We introduce the following function

$$
\omega_{\epsilon}(x)=
\cases
c_{\epsilon}e^{-\frac{\epsilon^{2}}{\epsilon^{2}-|x|^{2}}}
,& \hbox {if}~ |x|\leq \epsilon,\\
0 ,& \hbox{if}~ |x|>\epsilon.
\endcases
$$

  The constant is choosed $c_{\epsilon}$ such that the following equality
holds

$\int\omega_{\epsilon}(x)dx=1$.

Thus

$$
c_{\epsilon}\epsilon^{n}\int\limits_{|\xi|<1}e^{-\frac{1}{1-|\xi|^{2}}}d\xi=1.
$$

If $\chi(\cdot)$ be characteristic function of set
$K_{2\epsilon}$ then for any $\epsilon >0$ the  function

$$
\eta_{\epsilon}(x)=\int \chi(z)\omega_{\epsilon}(x-z)dz.
$$

satisfies the following relations [8, p.89]:

$$
0\leq\eta_{\epsilon}(x)\leq 1,\quad \eta_{\epsilon}(x)=1,~~  x\in K_{\epsilon},
$$
$$
\eta_{\epsilon}(x)=0,~~ x\notin K_{3\epsilon},\quad
\eta_{\epsilon}(x)\in C^{\infty}(R^{n}), \quad
 |\eta^{(\alpha)}_{\epsilon}(x)|
\leq L_{\alpha}\epsilon^{-|\alpha|}.\eqno(4)
$$

\bigskip

The next statement follows from the axiom of continuity:

\proclaim
{\bf Statement}
 If ~$\zeta$~ be random vector in space
~$R^{n}$ , then the following 

representation takes place

$$
 E\eta_{\epsilon}(\zeta)=P\left(\zeta\in K\right) + l_{\epsilon},
\quad \hbox{where}\quad l_{\epsilon}\geq 0,~~l_{\epsilon}\to 0,
~~\hbox{when}~~ \epsilon\to \infty.
\eqno(5)
$$
\endproclaim

\bigskip

\proclaim
{\bf Lemma 1}  If ~$P(\xi(0)\in K)=1$~
and for some number~$\epsilon_{0}>0$~and any numbers
~$0<\epsilon \leq \epsilon_{0}$~ the following inequality takes place

$$
E\eta_{\epsilon}(\xi(t))\geq  E\eta_{\epsilon}(\xi(0)),\quad t\geq 0.
$$

then the following equality is true

~~$P\left(\xi(t)\in K\right)=1,\quad t\geq 0.$

\endproclaim

\bigskip
\demo
{Proof} Let the condition of Lemma be fulfilled
but statement of Lemma don't fulfill. If statement of Lemma
don't fulfill then there exists such $t_{*}$ that

$$
P(\xi(t_{*})\in K)<1 \eqno(6).
$$

Futher according to the statement (5) and the condition of Lemma 1
we have the the following inequality in point ~$t_{*}$

$$
P(\xi(t_{*})\in K)+l_{3,\epsilon}\geq
P(\xi(0)\in K)+l_{2,\epsilon}.
$$

Letting  ~$\epsilon\to 0$~, we arrive at

$$
P(\xi(t_{*})\in K)\geq P(\xi(0)\in K)=1.
$$

The latter one contradicts to (6). This contradiction proves
the Lemma 1.
\enddemo

We make the following assumption:
the boundary of  ~$\partial K_{\epsilon}$ belongs to
class  $C^{l}$, when the following condition of smoothness
of boundary of $K_{\epsilon}$ holds  for ~$\epsilon<\epsilon_{0}$
under small $\epsilon_{0}>0$.
The intersection of boundary of set ~$K_{\epsilon}$  with ball
~$U(x,\epsilon),$

$~x\in K_{2\epsilon}:$

$$
 \triangle_{\epsilon}(x):=\bar U(x,\epsilon)\cap
\partial K_{\epsilon}
$$

is surface whose equation in local coordinates
~$(y_{1}, \dots, y_{n-1})$~ with origin of coordinates in point
~$x_{0}\in \triangle_{\epsilon}(x)$~ has form
~$y_{n}=\varphi(y_{1},\dots,y_{n-1})$.

The function  ~$\varphi$ belongs to class ~$C^{l}$~ in region
 ~$\bar D_{\epsilon}$~, which is projection of
~$\triangle_{\epsilon}(x)$

~ on the plane ~$y_{n}=0$.

Let us denote by $\nu(z)=(\nu_{i}(z), \quad i=1,2,...,n)$ the
unit vector of external normal to boundary ~$S$~ in point $z\in S $.

It is known, that if surface is given by relation
~$Q(y_{1},\dots,y_{n})=0$~, here ~$Q(\cdot)$~ 

be smooth function, then
the unit vector of normal~$\vec n$~ has the following form

$$
\vec n=(Q_{y_{i}}/\sqrt{\sum\limits_{k}Q_{y_{k}}^{2}}, ~~ i=\overline{1,n}).
$$

Thus, if the  ~$\varphi$~ is differentiable ,then the next reprezentation
for  $\nu(z)$~ takes place locally

$$
\nu(z)= \left(\frac{1}{\sqrt{1+\sum\limits_{i\leq n-1}\varphi_{y_{i}}^{2}}},
\frac{-\varphi_{y_{k}}}{\sqrt{1+\sum\limits_{i\leq n-1}\varphi_{y_{i}}^{2}}},
~~k=\overline{1,n-1}\right).
$$

Suppose now that the bondary ~$K_{\epsilon}, \epsilon_{0}\geq \epsilon\geq 0$
under some ~$\epsilon_{0}>0$~ belongs to class ~$C^{2}$.

\bigskip

\proclaim
{ Theorem 1} If the following conditions are fulfilled

1.~~ The functions $a(t,x), b_{k}(t,x),1\leq  k\leq n$, in addition to
properties  (2), under fixed ~$t$~
belong according to classes ~$C^{2}(R^{n})$,~$C^{3}(R^{n})$.

2. ~~ $\sup\limits_{s\geq 0}\sup\limits_{z\in S_{\epsilon}}
\left(\sum\limits_{i} b_{ji}(s,z)
\nu_{i}(z) \right)= o(\epsilon),~\epsilon\to 0.
\quad 1\leq j\leq n, \quad s\geq 0; $

3. ~~ $\overline{\lim\limits_{\epsilon\to 0}}\sup\limits_{s\geq 0}
\sup\limits_{z\in S_{\epsilon}}
\left(\sum\limits_{i}a_{i}(s,z)\nu_{i}(z)
  - \frac{1}{2}\sum\limits_{i,j,k}\frac{\partial b_{ki}(s,z)}
{\partial z_{j}}\nu_{i}(z) b_{kj}(s,z) \right)< 0.$

then  (3) takes place.

\endproclaim
\bigskip
\demo
{ Proof}
Applying the Ito's formula, we get the following equality
$$
E\eta_{\epsilon}(\xi(t)) - E\eta_{\epsilon}(\xi(0))=
E\int\limits_{0}^{t}A\eta_{\epsilon}(\xi(s))ds.
$$

Here

$$
 A :=\sum\limits_{i=1}^{n}a_{i}(s,x)\frac{\partial}{\partial x_{i}} +
\frac{1}{2}\sum\limits_{i,j = 1}^{n}\sigma_{ij}(s,x)\frac{\partial^{2}}
{\partial x_{i}\partial x_{j}}.
$$

The matrix $\sigma(s,x)=(\sigma_{ij}(x),\,1\leq i,j \leq n), $ is defined
in the following way

$$
\sigma(s,x)=B^{T}(s,x)B(s,x),\quad B(s,x):=(b_{ki}(s,x), 1\leq k,i\leq n).
$$

According to the Lemma 1 and definition of function ~$\eta_{\epsilon}(x)$~
for proof of invariance of set ~$K$~ it suffices to prove the
next inequality

$$
A\eta_{\epsilon}(x)\geq 0,\quad x\in
\left(K_{3\epsilon}\setminus K_{\epsilon}\right),
\quad s\geq 0.
$$

It is not difficult to check the following properties
of function  ~$\omega_{\epsilon}(x-z)$~

$$
-\frac{\partial}{\partial z_{i}}(\omega_{\epsilon} (x-z)) =
\frac{\partial}{\partial x_{i}}(\omega_{\epsilon} (x-z)); ~~
\frac{\partial^{2}}{\partial z_{i}\partial z_{j}}(\omega_{\epsilon}(x-z)) =
\frac{\partial^{2}}{\partial x_{i}\partial x_{j}}(\omega_{\epsilon}(x-z));
$$
$$
\omega_{\epsilon}(x-z)|_{|x-z|=\epsilon}
=\frac{\partial}{\partial x_{i}}(\omega_{\epsilon}(x-z))|_{|x-z|=\epsilon}=
\frac{\partial^{2}}{\partial x_{i}\partial x_{j}}
(\omega_{\epsilon}(x-z))|_{|x-z|=\epsilon} = 0.
$$

We define set

$$
K_{\epsilon}(x):=\left\{z: ~ |z-x| \leq \epsilon\cap
\left(K_{3\epsilon}\setminus
K_{\epsilon}\right) \right\}.
$$

Further, applying the properties of function  $\omega_{\epsilon}(x)$
and Taylor-series 

expansion of functions
~$a_{i}(s,x),~\sigma_{ij}(s,x)$~ in point ~~$z$~ we get

$$
A\eta_{\epsilon}(x) = \Bigl(\sum\limits_{i} a_{i}(s,x)\frac {\partial}
{\partial x_{i}}
+ \frac{1}{2}\sum\limits_{i,j}\sigma_{ij}(s,x)\frac{\partial^{2}}
{\partial x_{i}\partial x_{j}}\Bigr)\int\limits_{R^{n}}
\chi(z)\omega_{\epsilon}(x-z)dz =
$$
$$
= - \int\limits_{K_{\epsilon}(x)}\sum\limits_{i}\left(a_{i}(s,z)+
\sum\limits_{k}\frac{\partial a_{i}(s,z)}{\partial z_{k}}
(x_{k}-z_{k})+\right.
$$
$$
+\left.\frac{1}{2}\sum\limits_{k,j}\frac{\partial^{2}
a_{i}(s,\theta(x,z))}{\partial z_{k}\partial z_{j}}(x_{k}-z_{k})(x_{j}-z_{j})
\right)\frac{\partial}{\partial z_{i}}\omega_{\epsilon}(x-z)dz+
$$
$$
+ \int\limits_{K_{\epsilon}(x)}\frac{1}{2}
\sum\limits_{i,j}\left(\sigma_{ij}(s,z)+
\sum\limits_{k}
\frac{\partial\sigma_{ij}(s,z)}{\partial z_{k}}(x_{k}-z_{k}) +\right.
$$
$$
+\frac{1}{2}\sum\limits_{k,m}
\frac{\partial^{2}\sigma_{ij}(s,z)}
{\partial z_{k}\partial z_{m}}(x_{k}-z_{k})(x_{m}-z_{m})+
$$
$$
+\left. \frac{1}{6}\sum\limits_{k,m,l}\frac{\partial^{3}
\sigma_{ij}(s,\theta_{1}(x,z))}{\partial z_{k}\partial z_{m}\partial z_{l}}
(x_{k}-z_{k})(x_{m}-z_{m})(x_{l}-z_{l})\right)\frac{\partial^{2}
\omega_{\epsilon}(x-z)}{\partial z_{i}\partial z_{j}}dz=
$$

$$
= -
\int\limits_{K_{\epsilon}(x)}\sum\limits_{i}\Bigl\{\frac{\partial}{\partial
z_{i}}\left(a_{i}(s,z)\omega_{\epsilon}(x-z)\right)
-\omega_{\epsilon}(x-z)
\frac{\partial a_{i}(s,z)}{\partial z_{i}}+\Bigr.
$$

$$
+\sum\limits_{k}\frac{\partial}{\partial z_{i}}\left((x_{k}-z_{k})
\omega_{\epsilon}(x-z)\frac{\partial a_{i}(s,\theta(x,z))}
{\partial z_{k}}\right)
+\omega_{\epsilon}(x-z)\frac{a_{i}(s,z)}{\partial z_{i}}-
$$

$$
- \sum\limits_{k}(x_{k}-z_{k})\omega_{\epsilon}(x-z)\frac{\partial^{2}
 a_{i}(s,z)}{\partial z_{i}\partial z_{k}}+
$$
$$
+\Bigl.\frac{1}{2} \sum\limits_{k,j}\frac{\partial^{2}a_{i}(s,\theta(x,z))}
{\partial z_{k}\partial z_{j}}(x_{k}-z_{k})(x_{j}-z_{j})\frac{\partial
\omega_{\epsilon}(x-z)}{\partial z_{i}}\Bigr\}dz+
$$

$$
+\frac{1}{2}\int\limits_{K_{\epsilon}(x)}\sum\limits_{i,j}\Bigl\{
\frac{\partial}{\partial z_{i}}\left(\sigma_{ij}(s,z)
\frac{\partial \omega_{\epsilon}}{\partial z_{j}}\right)- \Bigr.
$$

$$
- \frac{\partial \sigma_{ij}(s,z)}{\partial z_{i}}
\frac{\partial\omega_{\epsilon}(x-z)}
{\partial z_{j}}+
\sum\limits_{k}\frac{\partial}{\partial z_{i}}\left(
(x_{k}-z_{k}) \frac{\partial \sigma_{ij}(s,z)}{\partial z_{k}}
\frac{\partial\omega_{\epsilon}(x-z)}{\partial z_{j}}\right)+
$$

$$
+ \frac{\partial \sigma_{ij}(s,z)}{\partial z_{i}}
\frac{\partial\omega_{\epsilon}(x-z)}{\partial z_{j}}-
\sum\limits_{k}(x_{k}-z_{k})\frac{\partial^{2}\sigma_{ij}(s,z)}{\partial z_{i}
\partial z_{k}}\frac{\partial\omega_{\epsilon}(x-z)}{\partial z_{j}}+
$$
$$
+\frac{1}{2}\sum\limits_{k,m}\frac{\partial}{\partial z_{i}}\left((x_{k}-z_{k})
(x_{m}-z_{m})
\frac{\partial^{2}
\sigma_{ij}(s,z)}{\partial z_{k}\partial z_{m}}
\frac{\partial \omega_{\epsilon}(x-z)}{\partial z_{j}}\right)+
$$

$$
+\frac{1}{2}\sum\limits_{m}(x_{m}-z_{m})
\frac{\partial^{2}
\sigma_{ij}(s,z)}{\partial z_{i}\partial z_{m}}
\frac{\partial \omega_{\epsilon}(x-z)}{\partial z_{j}}+
$$

$$
+\frac{1}{2}\sum\limits_{k}(x_{k}-z_{k})\frac{\partial^{2}
\sigma_{ij}(s,z)}{\partial z_{k}\partial z_{i}}
\frac{\partial \omega_{\epsilon}(x-z)}{\partial z_{j}}-
$$

$$
-\frac{1}{2}\sum\limits_{k,m}(x_{k}-z_{k})
(x_{m}-z_{m})
\frac{\partial^{3}
\sigma_{ij}(s,z)}{\partial z_{i}\partial z_{k}\partial z_{m}}
\frac{\partial \omega_{\epsilon}(x-z)}{\partial z_{j}}+
$$
$$
+\Bigl.\sum\limits_{k,m,l}\frac{\partial^{3}\sigma_{ij}(s,\theta_{1}(x,z))}
{\partial x_{k}\partial x_{m}\partial x_{l}}
(x_{k}-z_{k})(x_{m}-z_{m})(x_{l}-z_{l})\frac{\partial^{2}\omega_{\epsilon}(x-z)}
{\partial z_{i} \partial z_{j}}\Bigr\}dz.\eqno(7)
$$

Here ~$\theta(x,z)=z+\theta(x-z),~~\theta_{1}(x,z)=z+\theta_{1}(x-z);~~
 0\leq \theta_{1},\theta\leq 1$.

Applying Ostrogradskii - Gauss theorem we transform some integrals
in right part of  (7) to integrals on the surface ~$\triangle_{\epsilon}(x)$.
 Further, we obtain the estimate of smallness for some surface integrals
and the volume integrals.

We set

$$
f_{\epsilon}(i,x,z)=\frac{2\epsilon^{2}(x_{i}-z_{i})}{(\epsilon^{2}-
|x-z|^{2})^{2}},\quad\hbox{now}\quad
\frac{\partial \omega_{\epsilon}(x-z)}{\partial z_{i}}=
f_{\epsilon}(i,x,z) \omega_{\epsilon}(x-z).
$$

So

$$
\int\limits_{K_{\epsilon}(x)}\sum\limits_{i}
\frac{\partial}{\partial
z_{i}}\left(a_{i}(s,z)\omega_{\epsilon}(x-z)\right)dz=
 \int\limits_{\triangle_{\epsilon}(x)}
\sum\limits_{i}a_{i}(s,z)\nu_{i}(z)
\omega_{\epsilon}(x-z)
d\beta_{z}.\eqno(8)
$$

We make use of Cuachy- 'unyakovskii inequality in (9) and later on.

$$
|\int\limits_{K_{\epsilon}(x)}
\sum\limits_{i}
\frac{\partial}{\partial z_{i}}
\left(\sum\limits_{k}
(x_{k}-z_{k})
\omega_{\epsilon}(x-z)\frac{\partial a_{i}(s,z)}
{\partial z_{k}}\right)|=
$$
$$
=|\int\limits_{\triangle_{\epsilon}(x)}
\sum\limits_{i}
\nu_{i}(z)
\sum\limits_{k}
(x_{k}-z_{k})
\omega_{\epsilon}(x-z)
\frac{\partial a_{i}(s,z)}
{\partial z_{k}}
d\beta_{z}|\leq
$$

$$
\leq
 \int\limits_{\triangle_{\epsilon}(x)}
\sqrt{\sum\limits_{k}(x_{k}-z_{k})^{2}}\sqrt{\sum\limits_{k}
\left(\sum\limits_{i}\nu_{i}(z)\frac{\partial a_{i}(s,z)}{\partial z_{k}}
\right)^{2}}\omega_{\epsilon}(x-z)d\beta_{z}\leq \epsilon c_{1},
$$
$$
 ~~c_{1}<\infty.\eqno(9)
$$

$$
|\int\limits_{K_{\epsilon}(x)}
\sum\limits_{i}\sum\limits_{k}(x_{k}-z_{k})
\omega_{\epsilon}(x-z)
\frac{\partial^{2}
 a_{i}(s,z)}{\partial z_{i}\partial z_{k}}dz|\leq
$$
$$
\leq
\int\limits_{K_{\epsilon}(x)}
\sqrt{\sum\limits_{k}(x_{k}-z_{k})^{2}}
\sqrt{\sum\limits_{k}\left(\sum\limits_{i}
\frac{\partial^{2}
 a_{i}(s,z)}{\partial z_{i}\partial z_{k}}\right)^{2}}
\omega_{\epsilon}(x-z)dz\leq \epsilon c_{2},
$$
$$
c_{2}<\infty.\eqno(10)
$$

$$
|\int\limits_{K_{\epsilon}(x)}\sum\limits_{k,j}\frac{\partial^{2}a_{i}(\theta
(x,z))}{\partial z_{k}\partial z_{j}}(x_{k}-z_{k})(x_{j}-z_{j})
\frac{\partial \omega_{\epsilon}(x-z)}{\partial z_{i}}dz|\leq
$$
$$
\leq \int\limits_{K_{\epsilon}(x)}\sqrt{\sum\limits_{j,k,l}\left(\frac{\partial^{2}
a_{i}(s,\theta(x,z))}{\partial z_{k}\partial z_{j}}\right)^{2}}\sum\limits_{k}
\left(x_{k}-z_{k}\right)^{2}\times
$$
$$
\times\sqrt{\sum\limits_{i}f^{2}_{\epsilon}(i,x,z)}
\omega_{\epsilon}(x-z)dz\leq \epsilon c_{3},
\quad   c_{3}<\infty.\eqno(11)
$$

We exploit the condition (2) for estimate in (12)

$$
|\frac{1}{2}\int\limits_{K_{\epsilon}(x)}\sum\limits_{i,j}
\frac{\partial}{\partial z_{i}}\left(\sigma_{ij}(s,z)
\frac{\partial \omega_{\epsilon}}{\partial z_{j}}\right)dz|=
\frac{1}{2}|
 \int\limits_{\triangle_{\epsilon}(x)}
\sum\limits_{i}\nu_{i}(z)\sigma_{ij}(s,z)
\frac{\partial \omega_{\epsilon}}{\partial z_{j}}d\beta_{z}|=
$$

$$
=\frac{1}{2}|\int\limits_{\triangle_{\epsilon}(x)}\sum\limits_{k}
\sum\limits_{i}b_{ki}(s,z)\nu_{i}(z)
\sum\limits_{j}b_{kj}(s,z)
\frac{\partial \omega_{\epsilon}}{\partial z_{j}}d\beta_{z}|\leq
$$

$$
\leq \frac{1}{2}|
 \int\limits_{\triangle_{\epsilon}(x)}
|\sum\limits_{i}b_{ki}(s,z)\nu_{i}(z)|\sqrt{\sum\limits_{j}b^{2}_{kj}(s,z)}
\sqrt{\sum\limits_{j}f_{\epsilon}^{2}(j,x,z)}\times
$$
$$
\times\omega_{\epsilon}(x-z)d\beta_{z}=o(1).\eqno(12)
$$

We exploit the relation  $\sigma_{ij}(s,z)=\sigma_{ji}(s,z)$ for estimate
in (13)

$$
\frac{1}{2}|\int\limits_{K_{\epsilon}(x)}\sum\limits_{i,j}
\Bigl(-\sum\limits_{k}(x_{k}-z_{k})\frac{\partial^{2}\sigma_{ij}(s,z)}{\partial z_{i}
\partial z_{k}}\frac{\partial\omega_{\epsilon}(x-z)}{\partial z_{j}}+\Bigr.
$$

$$
+\frac{1}{2}\sum\limits_{m}(x_{m}-z_{m})
\frac{\partial^{2}
\sigma_{ij}(s,z)}{\partial z_{i}\partial z_{m}}
\frac{\partial \omega_{\epsilon}(x-z)}{\partial z_{j}}+
$$

$$
\Bigl.+\frac{1}{2}\sum\limits_{k}(x_{k}-z_{k})\frac{\partial^{2}
\sigma_{ij}(s,z)}{\partial z_{k}\partial z_{i}}
\frac{\partial \omega_{\epsilon}(x-z)}{\partial z_{j}}\Bigr)dz|=
$$
$$
=\frac{1}{2}|\int\limits_{\triangle_{\epsilon}(x)}\sum\limits_{i,j}\nu_{j}(z)
\Bigl(-\sum\limits_{k}(x_{k}-z_{k})\frac{\partial^{2}\sigma_{ij}(s,z)}{\partial z_{i}
\partial z_{k}}
+\frac{1}{2}\sum\limits_{m}(x_{m}-z_{m})
\frac{\partial^{2}
\sigma_{ij}(s,z)}{\partial z_{i}\partial z_{m}}+\Bigr.
$$

$$
\Bigl.+\frac{1}{2}\sum\limits_{k}(x_{k}-z_{k})
\frac{\partial^{2}
\sigma_{ij}(s,z)}{\partial z_{k}\partial z_{i}}\Bigr)\omega_{\epsilon}(x-z)
d\beta_{z}+
$$

$$
+\frac{1}{2}\int\limits_{K_{\epsilon}(x)}\sum\limits_{i,j}
\Bigl(
\frac{\partial^{2}\sigma_{ij}(s,z)}{\partial z_{i}
\partial z_{j}}-
\frac{1}{2}\frac{\partial^{2}\sigma_{ij}(s,z)}{\partial z_{i}
\partial z_{j}}-\frac{1}{2}
\frac{\partial^{2}
\sigma_{ij}(s,z)}{\partial z_{j}\partial z_{i}}\Bigr)\omega_{\epsilon}(x-z)
dz+
$$

$$
+\frac{1}{2}\int\limits_{K_{\epsilon}(x)}\sum\limits_{i,k,j}(x_{k}-z_{k})
\Bigl(-
\frac{\partial^{3}\sigma_{ij}(s,z)}{\partial z_{j}\partial z_{i}
\partial z_{k}}+
\frac{1}{2}\frac{\partial^{3}\sigma_{ij}(s,z)}{\partial z_{j}\partial z_{i}
\partial z_{k}}+ \frac{1}{2}\frac{\partial^{3}\sigma_{ij}(s,z)}
{\partial z_{j}\partial z_{k}\partial z_{i}}\Bigr)\omega_{\epsilon}(x-z)
dz|\leq
$$
$$
\leq\int\limits_{\triangle_{\epsilon}(x)}\sqrt{\sum\limits_{j}\nu_{j}^{2}}(z)
\sqrt{\sum\limits_{k}(x_{k}-z_{k})^{2}}\sqrt{\sum\limits_{i,j,k}
\Bigl(\frac{\partial^{2}\sigma_{ij}(s,z)}{\partial z_{i}\partial z_{k})}
\Bigr)^{2}}\omega_{\epsilon}(x-z)dz+
$$
$$
+\int\limits_{K_{\epsilon}(x)}
\sqrt{\sum\limits_{k}(x_{k}-z_{k})^{2}}\sqrt{\sum\limits_{i,j,k}
\Bigl(\frac{\partial^{3}\sigma_{ij}(s,z)}
{\partial z_{j}\partial z_{i}\partial z_{k})}
\Bigr)^{2}}
\omega_{\epsilon}(x-z)dz\leq c_{4} \epsilon,\quad c_{4}<\infty.\eqno(13)
$$

$$
\frac{1}{4}
|\int\limits_{K_{\epsilon}(x)}
\sum\limits_{i}\frac{\partial}{\partial z_{i}}\left(\sum\limits_{j}
\sum\limits_{k,m}(x_{k}-z_{k})
(x_{m}-z_{m})
\frac{\partial^{2}
\sigma_{ij}(s,z)}{\partial z_{k}\partial z_{m}}
\frac{\partial \omega_{\epsilon}(x-z)}{\partial z_{j}}\right)dz|=
$$

$$
=\frac{1}{4}|
\int\limits_{\triangle_{\epsilon}(x)}
\sum\limits_{i}\nu_{i}(z)\sum\limits_{j}
\sum\limits_{k,m}(x_{k}-z_{k})
(x_{m}-z_{m})
\frac{\partial^{2}
\sigma_{ij}(s,z)}{\partial z_{k}\partial z_{m}}
f_{\epsilon}(j,x,z) \times
$$

$$
\times\omega_{\epsilon}(x-z)d\beta_{z}|\leq
\frac{1}{4}
\int\limits_{\triangle_{\epsilon}(x)}
\sum\limits_{k}\left(x_{k}-z_{k}\right)^{2}
\sqrt{\sum\limits_{j}f_{\epsilon}^{2}(j,x,z)}
\sqrt{\sum\limits_{i}\nu^{2}_{i}(z)}\times
$$

$$
\times\sqrt{\sum\limits_{i,j,k,m}\left(
\frac{\partial^{2}
\sigma_{ij}(s,z)}{\partial z_{k}\partial z_{m}}\right)^{2}}
d\beta_{z}\leq \epsilon c_{5}, ~~ c_{5}<\infty.\eqno(14)
$$

$$
\frac{1}{4}|
\int\limits_{K_{\epsilon}(x)}
\sum\limits_{k,m}(x_{k}-z_{k})
(x_{m}-z_{m})
\frac{\partial^{3}
\sigma_{ij}(s,z)}{\partial z_{i}\partial z_{k}\partial z_{m}}
\frac{\partial \omega_{\epsilon}(x-z)}{\partial z_{j}}dz|\leq
$$
$$
\leq
\frac{1}{4}
\int\limits_{K_{\epsilon}(x)}
\sum\limits_{k}\left(x_{k}-z_{k}\right)^{2}\sqrt{\sum\limits_{j}
f_{\epsilon}^{2}(j,x,z)}\times
$$
$$
\sum\limits_{i}
\sqrt{\sum\limits_{j,k,m}
\left(\frac{\partial^{3}\sigma_{ij}(s,z)}
{\partial z_{i}\partial z_{k}\partial z_{m}}\right)^{2}}
\omega_{\epsilon}(x-z)dz\leq \epsilon c_{6},
~~c_{5}<\infty.\eqno(15)
$$

We exploit the following relation

$$
\frac{\partial^{2}\omega_{\epsilon}(x-z)}{\partial z_{i}\partial z_{j}}=
$$

$$
=\left(f_{\epsilon}(i,x,z)f_{\epsilon}(j,x,z) -\frac{2\epsilon^{2}\delta_{ij}}
{\left(\epsilon^{2}-|x-z|^{2}\right)^{2}}-\frac{8\epsilon^{2}(x_{i}-z_{i})
(x_{j}-z_{j})}{\left(\epsilon^{2}-|x-z|^{2}\right)^{3}}\right)
\omega_{\epsilon}(x-z),
$$
here $\delta_{ij}$ be Kronecker's symbol,

for estimate in  (16).

$$
\frac{1}{12}|\int\limits_{K_{\epsilon}(x)}\sum\limits_{i,j,k,m,l}
\frac{\partial^{3}\sigma_{ij}(s,\theta_{1}(x,z))}{\partial z_{k}
\partial z_{k}\partial z_{l}}(x_{k}-z_{k})(x_{m}-z_{m})(x_{l}-z_{l})
\frac{\partial^{2}\omega_{\epsilon}(x-z)}{\partial z_{i}\partial z_{j}}dz|\leq
$$
$$
\leq \frac{1}{12}\int\limits_{K_{\epsilon}(x)}\left(\sum\limits_{m}(x_{m}-z_{m})^{2}
\right)^{\frac{3}{2}}\sqrt{\sum\limits_{i,j,k,m,l}\left(\frac{\partial^{3}
\sigma_{ij}(s,\theta_{1}(x,z))}{\partial z_{k}\partial z_{m}\partial z_{l}}
\right)^{2}}\times
$$
$$
\times \left(\sum\limits_{i}f_{\epsilon}^{2}(i,x,z)+ \frac{2\epsilon^{2}\sqrt{n}}
{\left(\epsilon^{2}-|x-z|^{2}\right)^{2}}+\frac{8\epsilon^{2}\sum\limits_{i}
(x_{i}-z_{i})^{2}}{\left(\epsilon^{2}-|x-z|^{2}\right)^{3}}\right)
\omega_{\epsilon}(x-z)dz\leq \epsilon c_{6},
$$
$$
 c_{6}<\infty.\eqno(16)
$$

Concider remaining summand in (7).
$$
\frac{1}{2}
\int\limits_{K_{\epsilon}(x)}
\sum\limits_{k}\frac{\partial}{\partial z_{i}}\left(
(x_{k}-z_{k}) \frac{\partial \sigma_{ij}(s,z)}{\partial z_{k}}
\frac{\partial\omega_{\epsilon}(x-z)}{\partial z_{j}}\right)=
$$
$$
=\frac{1}{2}\int\limits_{\triangle_{\epsilon}(x)}
\sum\limits_{k}\sum\limits_{p}
\sum\limits_{j}\frac{\partial b_{pj}(z)}{\partial z_{k}}
\sum\limits_{i}\nu_{i}(z)b_{pi}(s,z)(x_{k}-z_{k})
\frac{\partial}{\partial z_{j}}\omega_{\epsilon}(x-z)d\beta_{z} +
$$

$$
+\frac{1}{2} \int\limits_{\triangle_{\epsilon}(x)}
\sum\limits_{k}\sum\limits_{p}
\sum\limits_{j}b_{pj}(s,z)
\sum\limits_{i}\nu_{i}(z)\frac{\partial b_{pi}(s,z)}
{\partial z_{k}}(x_{k}-z_{k})
\frac{\partial}{\partial z_{j}}\omega_{\epsilon}(x-z)d\beta_{z}=:
$$
$$
=: I_{1}+I_{2}.
$$

The first summand  $I_{1}$ is estimated analogy to (11)
with help condition 2 of theorem. Thus we have $I_{1}=o(1).$

We observe that by construction the points of boundary
$\partial \triangle_{\epsilon}(x)$ of set
~$\triangle_{\epsilon}(x)$~ 

have the following properties:

$$
z\in \partial \triangle_{\epsilon}(x)\Rightarrow |x-z|=\epsilon\Rightarrow
\omega_{\epsilon}(x-z)=0
$$

Now we make use of local property of surfase
~$\partial K_{\epsilon}$~ for more precise 

representation
of summand ~$I_{2}$~.

The variables $z_{i}, ~i=\overline{1,n}$~ in ~$\triangle_{\epsilon}(x)$~
 have form ~$z_{i}=y_{i},~i\leq n-1, z_{n}=
\varphi(y_{1},\dots,y_{n-1})$.

Put ~$\hat y=
(y_{1},\dots,y_{n}),~~\hbox{here}~ ~y_{n}=\varphi(y_{1},\dots,y_{n-1})$.

The domain  ~$\bar D_{\epsilon}(x)$~ which corresponds to
~$\triangle_{\epsilon}(x)$~has the following form

$$
\bar D_{\epsilon}(x)=\{(y_{1},\dots,y_{n-1}):
~~|x-\hat y|\leq \epsilon\}.\eqno(17)
$$

The boundary ~$\partial\bar D_{\epsilon}(x)$~ is set of points
for which in (17) the next equality is fulfilled.
Let ~$y'=(y_{1},\dots, y_{n-1})$.
Thus if  $y'\in \partial\bar D_{\epsilon}(x)$,
 then~$\omega_{\epsilon}(x-\hat y)=0.$

Set  ~$\omega_{\epsilon}(x-\hat y)=0$, under
$y'\notin \bar D_{\epsilon}(x)$.

Thus the function~$\omega_{\epsilon}(x-\hat y)$~
is finite function in space ~$R^{n-1}$~ with 

support
~$\bar D_{\epsilon}(x)$.
The following formula of integration by parts is true
for such functions  [8, p.106].

$$
f\in C^{1}\Rightarrow \int\limits_{\bar D_{\epsilon}(x)}
f\frac{\partial}{\partial y_{i}}
\omega_{\epsilon}(x-\hat y)dy'= -
\int\limits_{\bar D_{\epsilon}(x)}\omega_{\epsilon}(x-\hat y)
\frac{\partial}{\partial y_{i}}fdy'\quad i=\overline{1,n-1}.\eqno(18)
$$

Applying (18) to integration in ~$I_{2}$, we get

$$
2I_{2}=
$$
$$
= - \int\limits_{\bar D_{\epsilon}(x)}\sum\limits_{j}\sum\limits_{k}\frac{\partial}
{\partial y_{j}}\left(\sum\limits_{p,i}b_{pj}(s,\hat y)\nu_{i}(\hat y)
\frac{\partial}{\partial y_{k}}b_{pi}(s,\hat y)\right)(x_{k}-y_{k})
\omega_{\epsilon}(x-\hat y)dy'+
$$

$$
+\int\limits_{\bar D_{\epsilon}(x)}
\sum\limits_{j}\sum\limits_{k}
\sum\limits_{p,i}b_{pj}(s,\hat y)\nu_{i}(\hat y)\frac{\partial}
{\partial y_{k}}b_{pi}(s,\hat y)
\frac{\partial}{\partial y_{j}}y_{k}
\omega_{\epsilon}(x-\hat y)dy'=:
$$
$$
=: I_{21}+I_{22}.
$$

To estimate of summand ~$I_{21}$ with help 
Cauchy - Bunyakovskii's inequality
we will use the condition 1 and the 
supposition that surface belongs to class $C^{2}$

Later on it is convenient to omit the argument of functions.

$$
|I_{21}|\leq
\int\limits_{\bar D_{\epsilon}(x)}
\sum\limits_{j}\sum\limits_{p,i}
\left|\sum\limits_{k}
\frac{\partial}{\partial y_{j}}\left(b_{pj}\nu_{i} \frac{\partial}{\partial
 y_{k}}b_{pi}\right)(x_{k}-y_{k})\right|\omega_{\epsilon}(x-y)dy'\leq
$$
$$
\leq
\int\limits_{\bar D_{\epsilon}(x)}
\left(\sum\limits_{k}(x_{k}-y_{k})^{2}\right)^{\frac{1}{2}}
\sum\limits_{j}\sum\limits_{p,i}
\left(\sum\limits_{k}\left\{(\frac{\partial b_{pj}}{\partial y_{j}}\nu_{i}
\frac{\partial b_{pi}}{\partial y_{k}})^{2}+\right.\right.
$$
$$
\left.\left. +(b_{pj}\frac{\partial \nu_{i}}{\partial y_{j}}
\frac{\partial b_{pi}}
{\partial y_{k}})^{2}+ (b_{pj}\nu_{i}\frac{\partial^{2} b_{pi}}
{\partial y_{j}
\partial y_{k}})^{2}\right\}\right)^{\frac{1}{2}}\omega_{\epsilon}(x-\hat y)dy'\leq
\epsilon c_{21}.
$$

Here ~$c_{21}$~ is bounded constant.

Combining  (8)-(16) and latter one gives the following representation

$$
A\eta_{\epsilon}(x)=-\int\limits_{\bar D_{\epsilon}(x)}\left\{\sum\limits_{i}a_{i}(s,\hat y)
\nu_{i}(\hat y)-\right.
$$
$$
\left.-\frac{1}{2}\sum\limits_{j,k,p,i}b_{jp}(s,\hat y)\nu_{i}(\hat y)
\frac{\partial b_{pi}(s,\hat y)}{\partial y_{k}}\frac{\partial}{\partial y_{j}}
y_{k}\right\}\omega_{\epsilon}(x-\hat y)dy'+o(1).\eqno(19)
$$

Let $G(s,y)$ denote the second summ in braces of right part
of latter equality.

It takes place the following relation for new variables

$$
\frac{\partial}{\partial y_{n}}=
\sum\limits_{i=1}^{n-1}\varphi_{y_{i}}\frac{\partial}{\partial y_{i}}.
\eqno(20)
$$

It is not hard to calculate the following equalities for
partial derivatives in summands from ~$G(s,y)$~

$$
\frac{\partial b_{pi}}{\partial y_{k}}\frac{\partial}{\partial y_{j}}y_{k}=
\cases
\frac{\partial}{\partial y_{k}}b_{pi}
,& \hbox{if}~ j=k<n \\
0,&\hbox{if}\quad j\ne k ,~j<n,~k< n\\
\frac{\partial}{\partial y_{n}}b_{pn}\varphi_{x_{j}} ,& \hbox{if}~
j<n,~k=n\\
\frac{\partial b_{pi}}{\partial y_{k}}\varphi_{x_{k}},&
\hbox{if}~k<n,~j=n\\
\frac{\partial}{\partial y_{n}}b_{pi}\sum\limits_{m=1}^{n-1}\varphi_{y_{m}}^{2}
,& \hbox{if}~ j=k=n.
\endcases
$$

Now we will show that the function $G(s,y)$ coincides with the
following function from condition 3 of theorem completely

$$
  \sum\limits_{i,j,k}\frac{\partial b_{ki}(s,z)}
{\partial z_{j}}\nu_{i}(z)b_{kj}(s,z) , \quad s\geq 0;\eqno(21)
$$

in case when differentiation is fulfilled in coordinates ~$\hat y$.

Applying (20), we get the following equailities for differentiation in (21)

$$
\frac{\partial b_{pi}}{\partial z_{j}}=
\cases
\frac{\partial}{\partial y_{j}}b_{pi} +\frac{\partial b_{pi}}{\partial z_{n}}
\varphi_{x_{j}},,& \hbox{if}~ j<n \\
\sum\limits_{k}^{n-1}
\frac{\partial b_{ip}}{\partial y_{k}}\varphi_{x_{k}}+
\frac{\partial}{\partial z_{n}}\sum\limits_{l=1}^{n-1}\varphi_{x_{l}}^{2},&
\hbox{if}~ j=n.
\endcases
$$

It is clear that latter one defines the summands in  (21) which is identical
to the summands in the ~$G(s,y)$.

Thus it follows from representation (19) that under conditions theorem
there exists such  $\epsilon^{*}>0$ that the inequality
$A\eta_{\epsilon}(x)\geq 0$  is fulfilled for all
$\epsilon\leq \epsilon^{*}$. Theorem is prooved.
\enddemo

\bigskip

\head{References}
\endhead
\bigskip
\ref
\no 1
\by J.-P.Aubin and G.Da Prato
\paper Stochastic viability and invariance
\jour Ann. Scuola Norm. Sup. Pisa
\yr 1990
\vol l27
\pages 595-694
\endref
\ref
\no 2
\by R.Buckdahn , M.Quincampoix and A.Rascanu
\paper Viability property for a bacward stochastic
differential equation and application to partial
differential equations
\jour Probab.Theory Relat.Fields
\yr 2000
\vol 1l6
\page 485-504
\endref
\ref
\no 3
\by Il.I. Gikhman,  I.E.Klychkova
\paper Stochastic differential equations on the
embedded manifoldes
\jour Ukrain. math. journ.
\yr 1995
\vol 47
\page 174-179
\endref
\ref
\no 4
\by A.Milian
\paper Invariance for stohastic equations with regular coefficients
\jour Stochastic Analysis and Applications
\yr 1997
\vol 15
\page 91-101
\endref
\ref
\no 5
\by V.A.Gasanenko
\paper On invariant sets for stochastic differential equations
\jour Theory of stochastic processes
\yr 2003
\vol 9(25)
\page 60-64
\endref
\ref
\no 6
\by G.L.Kulinich, O.V. Pereguda
\paper The qualitative analysis of systems stochastic diffferential
Ito's equations
\jour Ukrain. math. journ.
\yr 2000
\vol 52
\page 1251-1256
\endref
\ref
\no 7
\by I.I. Gikhman, A.V.Skorokhod
\paper Inroduction to theory of random processes
\publ    Nauka
\publaddr Moskow
\yr 1977
\page 568 p.
\endref
\ref
\no 8
\by V.S.Vladimirov
\paper The equations of mathematical physics
\publ  Nauka
\publaddr Moskow
\yr 1988
\page 512 p.
\endref
\enddocument